\documentclass[11pt]{article}
\usepackage{amssymb}
\newtheorem{remark}{Remark}
\newtheorem{proposition}{Proposition}
\newtheorem{example}{Example}
\newcommand{\barF}{\overline{F}}
\newcommand{\barG}{\overline{G}}
\newcommand{\ds}{\displaystyle}
\def\LEQ#1{\mathop{{\leq}_{\hbox{\scriptsize\rm #1}}}}
\def\GEQ#1{\mathop{{\geq}_{\hbox{\scriptsize\rm #1}}}}
\title{\Large\bf A PARRONDO PARADOX IN RELIABILITY THEORY\footnote{Paper appeared 
in {\em The Mathematical Scientist} 32 (2007), no.\ 1, 17--22.}}
\author{
{\sc Antonio Di Crescenzo} \\
{\rm\small Dipartimento di Matematica e Informatica} \\
{\rm\small Universit\`a di Salerno} \\
{\rm\small 84084 Fisciano (SA), Italy} \\
{\rm\small E-mail: \tt adicrescenzo@unisa.it}
}
\date{\empty}
\begin{document}
\maketitle
\normalsize
\setlength{\baselineskip}{14pt}
\begin{abstract}
Parrondo's paradox arises in sequences of games in which a winning expectation
may be obtained by playing the games in a random order, even though each game in 
the sequence may be lost when played individually. We present a suitable version 
of Parrondo's paradox in reliability theory involving two systems in series, 
the units of the first system being less reliable than those of the second.
If the first system is modified so that the distributions of its new units 
are mixtures of the previous distributions with equal probabilities, then
under suitable conditions the new system is shown to be more reliable than
the second in the ``usual stochastic order'' sense.
 

\end{abstract}
%
%--------------------------------------------------------------------
\section{Introduction}
%--------------------------------------------------------------------
In their simpler formulation, Parrondo's games are composed of two `atomic' games: 
playing each atomic game individually the player is led to lose, whereas he may 
expect to win if the games are played randomly. These games were first described 
by Parrondo \cite{Pa96} at a workshop in 1996, and later studied by
Harmer and Abbott \cite{HaAb99a}, \cite{HaAb99b}, and Harmer {\em et al.\/}
\cite{HaAbTa00}. The counterintuitive behavior of Parrondo's games has
been exhibited under several sets of rules, depending on the player's current
capital or the history of the game. Parrondo's games can also be thought of 
as the discrete-time version of flashing Brownian ratchets, namely the
periodical or random combination of Brownian motion in an asymmetric potential
and of Brownian motion in a flat potential (see Parrondo {\em et al.\/}
\cite{PaHaAb00}, Astumian \cite{As97}, \cite{As01}, Allison and Abbott
\cite{AlAb01}, for instance). As with Brownian ratchets, Parrondo's games 
are gaining increasing attention in the literature: a feasible device for 
reproducing Parrondo's games has been proposed by Allison and
Abbott \cite{AlAb02} and the role of chaos in Parrondo's games has been
investigated by Arena {\em et al.\/} \cite{ArFaFoMa03}.  We also
refer the reader to the comprehensive review by Harmer and Abbott \cite{HaAb02},
and the list of articles shown in Web page \cite{WebPage}.
\par
The result implied by Parrondo's games, by which several losing strategies 
can be turned into a winning strategy by careful combination, is often referred 
to as ``Parrondo's paradox''. As one might expect, this surprising and intriguing 
fact has also kindled the interest of the popular press (see Blakeslee \cite{Bl00}). 
However, to avoid illusions it should be pointed out that Parrondo's phenomenon 
is not applicable to several real-life games (such as those played in casinos), 
or in the stock market. Indeed, as shown for instance by Philips and Feldman 
\cite{PhFe04}, Parrondo's paradox does not lead to any contradiction. The 
right combination of two losing strategies can be a winning strategy, without 
violating any mathematical of physical property. As in many other paradoxes 
in probability, the crucial intuition (or its lack) centers on the independence 
(or the lack of it) of the games. Parrondo's games are composed of 
{\em non-independent\/} losing games (which are very far from games played 
in the real world), and it should not be surprising that a suitable 
combination of non-independent losing games may produce a winning game. 
Nevertheless, even if Parrondo's paradox does not lead to contradictions, 
it does help to formulate and solve many interesting problems, such as 
finding improved or optimal strategies in sequences of games. 
\par
The aim of this note is to present a suitable version of Parrondo's paradox
in the field of reliability theory. We compare pairs of systems in series with 
two components, the units of the first system being less reliable than those of 
the second. We study the effect of a modification in the first system consisting 
in the random choice of its components: Assuming that each unit is chosen randomly 
from a set of components identical to the previous ones, then the distributions 
of the new units are mixtures of the previous distributions. 
With such randomization we obtain a system that,
under suitable conditions, is shown to be more reliable than the second one,
although the single starting units are less reliable than those of the second
system. This result, proved in Section 2, should not be considered contradictory, 
and indeed in the reliability theory literature it is well known
that mixtures of distributions do not preserve certain properties. 
For instance, whereas mixtures of decreasing failure rate distributions are always 
decreasing, very often mixtures of increasing failure rate distributions may  
also decrease (see Finkelstein and Esaulova \cite{FiEs01}, for instance).
\par
Finally, in Section 3 we give a sketch of a Parrondo's game based on the
outcomes of lifetimes of systems in series.
%--------------------------------------------------------------------
\section{The paradox}
%--------------------------------------------------------------------
Consider two systems, ${\cal S}_X$ and ${\cal S}_Y$, both formed by two
units connected in series. Let $X_1$ and $X_2$ be non-negative random
variables that describe the lifetimes of the units of system $S_X$.
We denote by $\barF_i(t)=\Pr\{X_i>t\}$, $t\geq 0$, the absolutely continuous
survival function of $X_i$, $i=1,2$. The random lifetime of ${\cal S}_X$ is
thus given by $X:=\min\{X_1,X_2\}$, and possesses the survival function
\begin{equation}
  \barF(t):=\barF_1(t)\,\barF_2(t), \qquad t\geq 0.
  \label{equation:1}
\end{equation}
Similarly, let $Y_1$ and $Y_2$ be non-negative random variables that describe
the lifetimes of the series units of system ${\cal S}_Y$. They are 
characterized
by absolutely continuous survival functions $\barG_1(t)$ and $\barG_2(t)$,
$t\geq 0$. Hence, the lifetime of ${\cal S}_Y$ is $Y:=\min\{Y_1,Y_2\}$, and
has the survival function
\begin{equation}
  \barG(t):=\barG_1(t)\,\barG_2(t), \qquad t\geq 0.
  \label{equation:8}
\end{equation}
The random variables $X_1,X_2,Y_1,Y_2$ are assumed to be independent, 
all with support
$(0,+\infty)$, i.e.\ with $\barF_i(t)>0$ and $\barG_i(t)>0$ for all 
$t\geq 0$, $i=1,2$.
\par
Furthermore, assume that the components of system ${\cal S}_Y$ are 
more reliable
than those of ${\cal S}_X$ with respect to the ``usual stochastic 
order'' (for a
detailed study of such stochastic order see Shaked and Shanthikumar 
\cite{ShSh94}):
we have $X_i\LEQ{st}Y_i$ for $i=1,2$, i.e.\ $\barF_i(t)\leq \barG_i(t)$ for 
all $t\geq 0$, $i=1,2$. In other words, for all times $t\geq 0$ lifetime 
$Y_i$ is more likely than $X_i$ to have a duration larger that $t$. 
From (\ref{equation:1}) and (\ref{equation:8}), we see that system ${\cal S}_Y$ 
is more reliable than ${\cal S}_X$, since $\barF(t)\leq \barG(t)$ for all 
$t\geq 0$, and then $X\LEQ{st}Y$.
\par
Let us now compare the above systems when the selection of the units 
of ${\cal S}_X$ is randomized. Assume that both components of the first 
system can be chosen randomly from two collections of units characterized 
by i.i.d.\ random lifetimes distributed as $X_1$ and $X_2$, respectively. 
In other words, we consider a new series system, ${\cal S}_{X^*}$ say, 
whose $i$-th component has survival function $\barF_i^*(t)$ given by 
a mixture with equal probabilities of the survival functions of $X_1$ 
and $X_2$, i.e.
$$
  \barF_i^*(t):={\barF_1(t)+\barF_2(t)\over 2},
  \qquad t\geq 0, \quad i=1,2.
$$
Consequently, the new system lifetime $X^*$ has survival function
\begin{equation}
  \barF^*(t):=\barF_1^*(t)\,\barF_2^*(t)
  =\left[{\barF_1(t)+\barF_2(t)\over 2}\right]^2,
  \qquad t\geq 0.
  \label{equation:2}
\end{equation}
 From (\ref{equation:1}) and (\ref{equation:2}) we have
$$
  \barF^*(t)-\barF(t)
  =\left[{\barF_1(t)-\barF_2(t)\over 2}\right]^2
  \geq 0,
  \qquad t\geq 0,
$$
so that $X^*\GEQ{st}X$. System ${\cal S}_{X^*}$ is thus more reliable 
than ${\cal S}_X$.
In other terms, the random choice of the components increases the system
lifetime in the usual stochastic order sense.
\par
Bearing in mind that lifetimes $X_1$ and $X_2$ are stochastically smaller 
than $Y_1$ and $Y_2$, respectively, we now face the following problem: to 
establish that, under certain conditions, system ${\cal S}_{X^*}$ is more 
reliable than ${\cal S}_Y$, i.e.:
\begin{equation}
\hbox{
Is it possible that $X^*\GEQ{st}Y$ even if $X_1\LEQ{st}Y_1$ and 
$X_2\LEQ{st}Y_2$?}
  \label{equation:7}
\end{equation}
In other words, recalling (\ref{equation:8}) and (\ref{equation:2}) 
we are looking
for survival functions $\barF_i(t)$ and $\barG_i(t)$, $i=1,2$, such that \\
{\em (i)} \ $[\barF_1(t)+\barF_2(t)]^2/4\geq \barG_1(t)\,\barG_2(t)$
for all $t\geq 0$, \\
{\em (ii)} \ $\barF_i(t)\leq \barG_i(t)$ for all $t\geq 0$, $i=1,2$.
\par
Hereafter we characterize the solutions of the above problem in
terms of a pair of functions that take values in the two-dimensional domain
belonging to the first quadrant of $\mathbb{R}^2$ and limited by the orthogonal
axis and by a time-varying branch of an equilateral hyperbola. 
Indeed, by setting
$$
  x_i(t)={\barG_i(t)\over \barF_i(t)}-1, \qquad t\geq 0,\quad i=1,2
$$
relations {\em (i)} and {\em (ii)} are satisfied if and only if,
for any fixed $t\geq 0$, we have
\begin{equation}
  x_1(t)\geq 0, \qquad
  x_2(t)\geq 0, \qquad
  x_1(t)+x_2(t)+x_1(t)\,x_2(t)\leq A(t),
  \label{equation:3}
\end{equation}
where we have set
$$
  A(t):={[\barF_1(t)-\barF_2(t)]^2\over 4\,\barF_1(t)\,\barF_2(t)},
  \qquad t\geq 0.
$$
\par
The answer to question (\ref{equation:7}) is affirmative. Indeed, two examples
of non-empty families of solutions are given below.
%
% =====================  proposition 1  =====================
\begin{proposition}
Two sets of sufficient conditions such that inequalities
(\ref{equation:3}) hold are the following: \\
{\em (a)} \ $x_1(t)=0$ and $0\leq x_2(t)\leq A(t)$, i.e.
\begin{equation}
  \barG_1(t)=\barF_1(t)
  \qquad \hbox{and} \qquad
  \barF_2(t)\leq \barG_2(t)
  \leq {[\barF_1(t)+\barF_2(t)]^2\over 4\,\barF_1(t)};
  \label{equation:5}
\end{equation}
{\em (b)} \ $x_1(t)=x_2(t)$ and $0\leq x_2(t)\leq \sqrt{A(t)+1}-1$, i.e.
\begin{equation}
  {\barG_1(t)\over\barF_1(t)}={\barG_2(t)\over\barF_2(t)}
  \qquad \hbox{and} \qquad
  \barF_1(t)\leq \barG_1(t)
  \leq {\barF_1(t)+\barF_2(t)\over 2}\,\sqrt{\barF_1(t)\over\barF_2(t)}.
  \label{equation:6}
\end{equation}
\end{proposition}
\par
Note that the functions appearing in the right-hand upper bounds of
(\ref{equation:5}) and (\ref{equation:6}) are not necessarily 
survival functions.
%
% =====================  REMARK 1  =====================
\begin{remark}{\rm
Denoting the probability density functions of $X_i$ and $Y_i$ respectively 
by $f_i(x)$ and $g_i(x)$, $i=1,2$, we give two necessary conditions such that 
inequalities (\ref{equation:3}) hold:
\begin{equation}
  f_1(0)=g_1(0)
  \qquad \hbox{and} \qquad
  f_2(0)=g_2(0).
  \label{equation:4}
\end{equation}
Indeed, by setting $h(t):=[\barF_1(t)+\barF_2(t)]^2/4-\barG_1(t)\,\barG_2(t)$
and $k_i(t):=\barG_i(t)-\barF_i(t)$, $i=1,2$, conditions {\em (i)}
and {\em (ii)} imply the non-negativity of $h(t)$ and $k_i(t)$, $i=1,2$;
moreover, since $h(0)=k_1(0)=k_2(0)=0$, by imposing the non-negativity
as $t\to 0^+$ of the derivatives of $h(t)$ and $k_i(t)$, $i=1,2$, we
immediately obtain (\ref{equation:4}).
}\end{remark}
%
% =====================  ESEMPIO 1  =====================
%
\begin{example}{\rm
A case in which conditions (\ref{equation:5}) are satisfied is
the following:
$$
  \barF_1(t)=\barG_1(t)=e^{-\lambda t},
  \qquad
  \barF_2(t)=(1+\lambda t)\,e^{-\lambda t},
  \qquad
  \barG_2(t)=\left[\left(1+{\lambda t\over 2}\right)^2
  -\left({\nu t\over 2}\right)^2\right]e^{-\lambda t},
$$
for $t\geq 0$, where $\lambda>0$ and $0\leq \nu\leq \lambda$. In this case
conditions (\ref{equation:4}) are satisfied, since $f_1(0)=g_1(0)=\lambda$
and $f_2(0)=g_2(0)=0$. The mean values of the component lifetimes are given by
$$
  E(X_1)=E(Y_1)={1\over \lambda},
  \qquad
  E(X_2)={2\over \lambda},
  \qquad
  E(Y_2)={2\over \lambda}+{1\over 2\lambda}\left[1-\left({\nu\over 
\lambda}\right)^2\right],
$$
so that $E(X_i)\leq E(Y_i)$, $i=1,2$. Moreover, we have that 
\begin{equation}
  E(X)={3\over 4\lambda},
  \qquad
  E(Y)={13\over 16\lambda}-{1\over 16\lambda}\left({\nu\over \lambda}\right)^2,
  \qquad
  E(X^*)={13\over 16\lambda}.
  \label{equation:10}
\end{equation}
Hence, we obtain $E(X)\leq E(Y)\leq E(X^*)$ for all $0\leq \nu\leq 
\lambda$, with
$E(X)<E(Y)=E(X^*)$ when $\nu=0$ and $E(X)=E(Y)<E(X^*)$ when $\nu=\lambda$.
}\end{example}
%
% =====================  ESEMPIO 2  =====================
%
\begin{example}{\rm
Let $\lambda>0$, and let $u(t)$ be a non-negative function for $t\geq 
0$, with $u(0)=1$
and such that ${{\rm d}\over {\rm d}t}u(t)\leq \lambda\,u(t)/2$, for $t\geq 0$.
By taking survival functions that satisfy the following relations:
$$
  \barF_1(t)=\left[u(t)\right]^2 e^{-\lambda t},
  \qquad
  \barF_2(t)=e^{-\lambda t},
$$
$$
  \left[u(t)\right]^2 e^{-\lambda t}\leq \barG_1(t)
  \leq u(t)\,{1+\left[u(t)\right]^2\over 2}\,e^{-\lambda t},
  \qquad
  \barG_2(t)={\barG_1(t)\over \left[u(t)\right]^2},
$$
for $t\geq 0$, conditions (\ref{equation:6}) are fulfilled. For instance,
by setting $u(t)=1+\lambda t/2$ and suitably choosing $\barG_1(t)$ we have:
\begin{equation}
  \begin{array}{cc}
  \barF_1(t)=\left(1+\ds{\lambda\,t\over 2}\right)^2\,e^{-\lambda t},
  &
  \barF_2(t)=e^{-\lambda t}, \\
  {} & {} \\
  \barG_1(t)=\left(1+\ds{\lambda t\over 2}\right)
  \left[1+\ds{\lambda t\over 2}+{(\nu t)^2\over 8}\right]e^{-\lambda t},
  & \quad
  \barG_2(t)=\left[1+\ds{(\nu t)^2/8\over 1+(\lambda t)/2}\right]e^{-\lambda t},
  \end{array}
  \label{equation:9}
\end{equation}
for $t\geq 0$, where $\lambda>0$ and $0\leq \nu\leq \lambda$.
Furthermore, from survival functions (\ref{equation:9}) it follows that
$f_1(0)=g_1(0)=0$ and $f_2(0)=g_2(0)=\lambda$, in agreement with 
(\ref{equation:4}).
}\end{example}
%
% =====================  REMARK 2  =====================
\begin{remark}{\rm
The identity at the left-hand-side of (\ref{equation:6}) yields a relation 
between the hazard rate functions of the component random lifetimes. Indeed, 
denoting by
$$
  h_{\barF}(t)=-{{\rm d}\over {\rm d}t}\ln \barF(t), \qquad t\geq 0
$$
the hazard rate function of a random lifetime characterized by the 
absolutely continuous survival function $\barF(t)$, we have
$$
  {\barG_1(t)\over\barF_1(t)}={\barG_2(t)\over\barF_2(t)}
  \qquad \hbox{for all $t\geq 0$}
$$
if and only if
$$
  h_{\barG_1}(t) - h_{\barF_1}(t) = h_{\barG_2}(t) - h_{\barF_2}(t)
  \qquad \hbox{a.e.}
$$
This relation is clearly satisfied by the hazard rates corresponding to
survival functions (\ref{equation:9}).
}\end{remark}
\par
We conclude this section by pointing out again that the affirmative answer
to question (\ref{equation:7}) is not a paradox. Indeed, assumptions
$X_1\LEQ{st}Y_1$ and $X_2\LEQ{st}Y_2$ imply $X^*\LEQ{st}Y^*$, and this
does not exclude the possibility that $X^*\GEQ{st}Y$ in some cases.
%--------------------------------------------------------------------
\section{A Parrondo's game}
%--------------------------------------------------------------------
On the basis of the above results, and assuming more realistically that the
independent random variables $X_1,X_2,Y_1,Y_2$ have finite mean values, we can
build up a game in which a player gains the outcome of $X-Y$. Hence, he wins
if a realization of $X$ is larger than a realization of $Y$. Unfortunately,
assumptions $X_1\LEQ{st}Y_1$ and $X_2\LEQ{st}Y_2$ lead to $X\LEQ{st}Y$, and
thus to a non-positive player's expected gain $E(X-Y)\leq 0$. We recall 
that in this case the allocation of units of system ${\cal S}_X$ is fixed
deterministically. However, let us assume that the game rules allow us to 
select randomly the units of ${\cal S}_X$, by means of mixtures with equal 
probabilities. This produces a new system lifetime $X^*$ with survival function 
given in (\ref{equation:2}), so that the player now gains the outcome of $X^*-Y$.
If the survival functions of $X_1,X_2,Y_1,Y_2$ satisfy conditions {\em (i)}
and {\em (ii)} then $X^*\GEQ{st}Y$, and the player's expected gain is
non-negative, being $E(X^*-Y)\geq 0$. For instance, with reference to the
case treated in Example 1, assuming $\lambda=1$ and $\nu=1/2$ from
(\ref{equation:10}) we have
$$
  E(X)={3\over 4},
  \qquad
  E(Y)={51\over 64},
  \qquad
  E(X^*)={13\over 16}.
$$
Hence, the player's expected gain is
$$
  E(X-Y)=-{3\over 64} \quad
  \hbox{when the units of the first system are deterministically set}
$$
and
$$
  E(X^*-Y)={1\over 64} \quad
  \hbox{when the units of the first system are randomly mixed}.
$$
Finally, this is a further confirmation of the standard conclusion of
Parrondo's games: even if fixed settings lead to losses, we may expect 
a win in the presence of random mixing.
%
%--------------------------------------------------------------------
\subsection*{Acknowledgement}
%--------------------------------------------------------------------
The author wishes to thank an anonymous referee for his careful reading of the paper. 
%--------------------------------------------------------------------

%

\end{document}